\documentclass[12pt,reqno]{amsart}
\usepackage{amscd}
\usepackage{amssymb}
\usepackage[mathscr]{eucal}
\usepackage{enumerate}

\makeatletter

\def\subsection{\@startsection{subsection}{2}%
  \z@{.5\linespacing\@plus.7\linespacing}{1pt}%
      {\normalfont\bfseries}}

\def\l@section{\@tocline{1}{0pt}{1pc}{}{\bfseries}}
\def\l@subsection{\@tocline{2}{0pt}{3.1em}{5pc}{}}

\makeatother

\allowdisplaybreaks[1]

\newtheorem{thm}{Theorem}[section]
\newtheorem{lem}[thm]{Lemma}
\newtheorem{cor}[thm]{Corollary}

\newtheorem{thm*}{Theorem}

\theoremstyle{definition}

\newtheorem{rem}[thm]{Remark}
\newtheorem{defn}[thm]{Definition}

\newtheorem{ex}[thm]{Example}

\theoremstyle{remark}

\numberwithin{equation}{section}

\newcommand{\C}{{\mathbb{C}}}
\newcommand{\R}{{\mathbb{R}}}

\newcommand{\N}{{\mathbb{N}}}

\DeclareMathOperator{\id}{id}

\DeclareMathOperator{\spa}{span}

\DeclareMathOperator{\Irr}{Irr}

\DeclareMathOperator{\Hom}{Hom}              % homs, endos, autos etc. 

\DeclareMathOperator{\End}{End}
\DeclareMathOperator{\Aut}{Aut}
\DeclareMathOperator{\Sect}{Sect}

\def\mP{\mathcal{P}}

\def\meA{\mathscr{A}}

\def\meH{\mathscr{H}}

\def\meK{\mathscr{K}}
\def\meL{\mathscr{L}}
\def\meM{\mathscr{M}}

\def\al{\alpha}
\def\be{\beta}
\def\ga{\gamma}
\def\de{\delta}

\def\vep{\varepsilon}
\def\ph{{\phi}}
\def\ps{{\psi}}
\def\vph{\varphi}

\def\om{\omega}
\def\si{\sigma}
\def\ta{\tau}
\def\th{\theta}

\def\Ps{\Psi}

\def\La{\Lambda}

\def\el{\ell}

\def\ovl{\overline}                 % overline, widehat, widetilde, opposite
\def\wdh{\widehat}
\def\wdt{\widetilde}

\def\subs{\subset}                   % for sets 

\def\setm{\setminus}
\def\nin{\notin}

\def\oti{\otimes}                    % tensor products
\def\col{\colon}
\def\ra{\hspace{-0.5mm}\rightarrow\!}

\def\bG{\mathbb{G}}
\def\bH{\mathbb{H}}

\def\ltG{L^2(\bG)}
\def\AG{A(\bG)}
\def\CG{C(\bG)}

\def\lG{L^\infty(\bG)}                      % function algebra on G 
\def\IG{\Irr(\bG)}

\def\suq{SU_q(2)}

\begin{document}

\title{A Galois correspondence for compact quantum group actions}

\author{Reiji Tomatsu}
\address{Department of Mathematical Sciences,
University of Tokyo, Komaba, \mbox{153-8914}\\
\indent JAPAN,}
\address{Department of Mathematics, 
K.U. Leuven, Celestijnenlaan 200B,
B-3001\\ 
\indent BELGIUM}

\email{tomatsu@ms.u-tokyo.ac.jp}

\thanks{Supported by JSPS}


\begin{abstract}
We establish a Galois correspondence for
a minimal action of a compact quantum group
${\mathbb G}$
on a von Neumann factor $M$.
This extends the result of Izumi, Longo and Popa
who treated the case of a Kac algebra.
Namely, there exists a one-to-one correspondence
between the lattice of left coideals of ${\mathbb G}$
and that of intermediate
subfactors of $M^{\mathbb G}\subset M$.
\end{abstract}

\maketitle

\section{Introduction}

The main purpose of this paper 
is to present a Galois correspondence for compact quantum group actions. 
The theory of Galois correspondences for group actions on von Neumann algebras 
was initiated by 
M. Nakamura and Z. Takeda \cite{NT1,NT2} and has been studied extensively 
in various settings 
by many researchers. 
Roughly speaking, 
the Galois correspondence for a group action $G$
on a von Neumann algebra $M$
refers to a one-to-one correspondence 
between the subgroups $H$ of $G$ and
the intermediate von Neumann subalgebras $M^G \subset M^H \subset M$, 
where $M^H$ denotes the fixed point algebra by the $H$-action. 

In \cite{NT2}, the Galois correspondence was established 
for a minimal action of a finite group on a II$_1$ factor. 
For compact group actions,  
A. Kishimoto obtained a Galois correspondence between 
normal closed subgroups and 
intermediate von Neumann subalgebras 
that are globally invariant under 
the compact group, 
assuming a certain condition on actions 
which is satisfied by minimal ones \cite{Ki}.
 
Another kind of Galois correspondence was provided 
in \cite{NT1} for 
crossed product von Neumann algebras 
in the case of a discrete group $G$ acting freely on a finite factor $M$. 
Their result again shows a one-to-one correspondence between 
the lattice of subgroups and that of intermediate subfactors
of $M \subset M \rtimes G$.
In \cite{Ta2}, M. Takesaki studied
a generalization of this result for
a locally compact abelian group action.
Y. Nakagami strengthened the result to the case of 
general locally compact group actions \cite{Na}.
In \cite{Ch}, 
H. Choda investigated 
the crossed product by free actions of discrete groups 
on a factor of arbitrary type 
and obtained the Galois correspondence for intermediate 
von Neumann subalgebras which are the ranges of
normal conditional expectations. 

M. Izumi, R. Longo and S. Popa \cite{ILP} have further developed theory 
of Galois correspondence to compact group minimal actions 
and 
discrete group free actions on factors of arbitrary type. 
Moreover, they obtained the Galois correspondence
for minimal actions of compact Kac algebras,
which unifies 
the results for compact groups and discrete groups. 

Therefore, it is natural 
to explore a generalization of the Galois correspondence of \cite{ILP}
to minimal actions of compact quantum groups. 
In \cite{E}, M. Enock focused on this problem in a more general setting, 
but there is a flaw in his proof. 
We point out the reason why the proof of \cite{ILP} 
does not work in the quantum case. 
The main step of their proof is 
to show the existence of a normal conditional expectation 
from an ambient von Neumann algebra onto any intermediate subfactor. 
However, this is no longer true in general for compact quantum group actions. 
Indeed, we consider a minimal action of $SU_q(2)$ ($0<q<1$) on a factor, 
whose existence has been shown by Y. Ueda \cite{U1}. 
The intermediate subfactor associated with 
a Podle\'{s} sphere $S_{q,\th}^2$ \cite{Po1,Po2}  
is not the range of a normal conditional expectation. 

For our purpose, we proceed to study an irreducible 
inclusion of discrete type in the sense of \cite{ILP}.
Let $N\subs M$ be an irreducible inclusion of discrete type. 
In \cite{ILP}, assuming a technical condition on the modular automorphism 
groups, it is shown that an intermediate subfactor $N\subs L\subs M$ 
is generated by $N$ and Hilbert spaces in $L$ which implement 
irreducible endomorphisms on $N$. 
We will strengthen the result 
by dropping assumptions on the modular automorphism groups. 
Our main result is the following (Theorem \ref{thm: L-discrete}): 

\begin{thm*}
Let $N\subs M$ be an irreducible inclusion of discrete type 
with the faithful normal conditional expectation $E\col M\ra N$. 
Let $N\subs L\subs M$ be an intermediate subfactor. 
Then the following statements hold.
\begin{enumerate}
\item 
The subfactor $N\subs L$ is discrete. 

\item 
Suppose $N$ is infinite. 
Let $\ga_N^M\col M\ra N$ and $\ga_N^L\col L\ra N$ be 
the canonical endomorphisms 
for $N\subs M$ and $N\subs L$, respectively. 
Then $[\ga_N^M|_N]$ contains $[\ga_N^L|_N]$ in $\Sect(N)$. 
\end{enumerate} 
\end{thm*}

The second statement of the above theorem is equivalent to saying 
that the bimodule ${}_N L^2(M)_N$ contains ${}_N L^2(L)_N$. 
This statement might sound trivial as in the case of finite $M$, 
but indeed we need a little more efforts to prove it 
(see Remark \ref{rem: bimodule}). 

Next we apply this result to minimal actions of compact quantum groups 
and prove the following Galois correspondence (Theorem \ref{thm: Galois}) 
which generalizes the correspondence presented in \cite{ILP}. 

\begin{thm*}
Let $\bG$ be a compact quantum group and $M$ a factor. 
Let $\al\col M\ra M\oti\lG$ be a minimal action. 
Then there exists a one-to-one correspondence between the lattice of 
left coideals in $\lG$ and that of 
intermediate subfactors of $M^\al\subs M$. 
\end{thm*}

\section{Irreducible inclusions of discrete type}

First, we fix notations. 
Throughout this paper, we always assume separability of von Neumann 
algebras. 
Let $M$ be a von Neumann algebra with predual $M_*$. 
For a weight $\ph$ on $M$, 
we set $n_\ph=\{x\in M\mid \ph(x^* x)<\infty\}$. 
The GNS representation of $M$ with respect to a faithful 
normal semifinite weight $\ph$ is denoted 
by the pair $\{H_\ph,\La_\ph\}$, where 
$\La_\ph\col n_\ph\ra H_\ph$ is the canonical injection 
to the GNS Hilbert space $H_\ph$. 
We always regard $M$ as a von Neumann subalgebra in $B(H_\ph)$. 

Let $N\subs M$ be an inclusion of von Neumann algebras. 
We denote by $\mP(M,N)$ 
the set of faithful normal semifinite operator valued weights 
from $M$ to $N$. 
For theory of operator valued weights,
readers are referred to \cite{H1,H2}.

For a subset $X\subs M$, we denote by $\ovl{X}^{\rm{w}}$,
$\ovl{{\rm co}}^{{\rm w}}(X)$ and $\ovl{\spa}^{\rm{w}}(X)$
the weak closure of $X$, the weak closure of the convex hull of $X$
and the weak closure of the linear space spanned by $X$,
respectively.

We denote by $\oti$ the minimal tensor product for $C^*$-algebras 
and the spatial tensor product for von Neumann algebras. 

We say that $\meH$ is a Hilbert space in $M$ if 
$\meH\subs M$ is a $\si$-weakly closed subspace 
and $\eta^* \xi\in\C$ for all $\xi,\eta\in\meH$ 
\cite{R}. 
The smallest projection $e\in M$ with $e\meH=\meH$ is 
called the support of $\meH$. 

We denote by $\End(M)$ and $\Sect(M)$ the set of 
endomorphisms and sectors on $M$, that is, $\Sect(M)$ 
is the set of equivalence classes of endomorphisms on $M$ 
by unitary equivalence. 
For the sector theory, readers are referred to \cite{I,L1,L2}. 

\subsection{Irreducible inclusions of discrete type} 

In this subsection, we recall the notion of discreteness
introduced in \cite{ILP}
for an inclusion of factors
and summarize the basic properties. 

Let $N\subs M$ be an inclusion of factors 
with a faithful normal conditional expectation $E_N^M\col M\ra N$. 
Fix a faithful state $\om\in N_*$. 
We set $\vph:=\om\circ E_N^M\in M_*$. 
Let $\{H_\vph,\La_\vph\}$ be the GNS representation 
of $M$ associated with the state $\vph$. 
We define the Jones projection $e_N\in B(H_\vph)$ by 
$e_N \La_\vph(x)=\La_\vph(E_N^M(x))$ for $x\in M$ 
and set the basic extension $M_1:=M\vee \{e_N\}''\subs B(H_\vph)$. 
The dual operator valued weight of $E_N^M$ is 
denoted by $\wdh{E}_M^{M_1}$ which is an element of $\mP(M_1,M)$
\cite{Ko1}. 
By definition, we have
$\wdh{E}_M^{M_1}(a e_N b)=ab$ for $a,b\in M$.
Define the faithful normal semifinite weight 
$\vph_1:=\vph\circ \wdh{E}_M^{M_1}$ on $M_1$.

In \cite[Definition 3.7]{ILP}, 
the discreteness of an inclusion of factors is introduced as follows. 
\begin{defn}
An inclusion of factors $N\subs M$ is said to be \emph{discrete}
when there exists a faithful normal conditional expectation
$E_N^M\col M\ra N$ such that its dual operator valued weight
$\wdh{E}_M^{M_1}$ is semifinite
on $N'\cap M_1$.
\end{defn}

Note that the discreteness is equivalent to saying that 
the $N$-$N$-bimodule $L^2(M)$ is the direct sum of 
$N$-$N$-bimodules of finite index. 

Let $N\subs M$ be an irreducible inclusion of discrete type. 
Then $N'\cap M_1$ can be decomposed into a direct sum 
of matrix algebras as \cite[Proposition 2.8]{ILP}: 
\[
N'\cap M_1=\bigoplus_{\xi\in\Xi}A_\xi, 
\]
where $A_\xi$ is a type I$_{\,n_\xi}$ factor for some $n_\xi\in\N$. 

Consider the case that $N$ is infinite. 
Let $\ga_N^M\col M\ra N$ be the canonical endomorphism 
of the inclusion $N\subs M$ \cite{L1,L2,L3}. 
By definition, we have the isomorphism of $M$-$N$-bimodules 
\[
{}_M{}_{\ga_N^M} L^2(N)_N \cong {}_M L^2(M)_N,
\]
where $L^2(M)$ and $L^2(N)$ denote 
the standard Hilbert space for $M$ and $N$, respectively. 
Since $N'\cap M_1=\End({}_N L^2(M)_N)$, 
a minimal projection $e_\xi$ in $A_\xi$ 
corresponds to an irreducible endomorphism with finite index, 
$\rho_\xi\in\End(N)$, that is, 
we have the following isomorphism of $N$-$N$-bimodules 
\cite[Lemma 3.1]{ILP}: 
\[
{}_N{}_{\rho_\xi} L^2(N)_N\cong {}_N\, e_\xi L^2(M)_N. 
\]
By definition, the sector $[\ga_N^M|_N]$ contains the sector $[\rho_\xi]$ 
with multiplicity $n_\xi$ in $\Sect(N)$. 

We define the subspace $\meH_\xi\subs M$ by 
\[
\meH_\xi=\{V\in M\mid Vx=\rho_\xi(x)V\ \mbox{for all}\ x\in N\}. 
\]
Then by the inner product $\langle V|W\rangle1=W^*V$, $\meH_\xi$ is 
a Hilbert space in $M$. 
We note that the support projection of $\meH_\xi$ 
may not be equal to $1$. 
We prepare another inner product defined by 
$(V,W)1=d(\xi)E_N^M(VW^*) $ for $V,W\in\meH_\xi$, 
where $d(\xi)$ is the square root of the minimal index 
of $\rho_\xi(N)\subs N$ \cite{Hi}. 
By \cite[Theorem 3.3]{ILP}, 
we have $\dim \meH_\xi=n_\xi$ and 
$A_\xi=\meH_\xi^* e_N \meH_\xi$. 

\subsection{Intermediate subfactors} 
Let $N\subs M$ be an irreducible inclusion of discrete type 
with infinite $N$.
Let $N\subs L\subs M$ be an intermediate subfactor.
We denote by $E_N^L$ the restriction of $E_N^M$ on $L$. 
For $\xi\in\Xi$, we define the Hilbert space $\meK_\xi$ in $L$ by 
\[
\meK_\xi:=\meH_\xi\cap L. 
\]
We set $\Xi_L :=\{\xi\in\Xi\mid \meK_\xi\neq0\}$ and 
$m_\xi:=\dim(\meK_\xi)$ for $\xi\in\Xi_L$. 

For $\xi\in\Xi$, 
we take a basis $\{V_{\xi_i}\}_{i\in I_\xi}$ in $\meH_\xi$ such that
$(V_{\xi_i},V_{\xi_j})=d(\xi)\de_{i,j}$.
If $\xi\in\Xi_L$, we may assume that the family
$\{V_{\xi_i}\}_{i\in I_\xi}$ contains a basis of $\meK_\xi$, 
which we denote by $\{V_{\xi_i}\}_{i\in I_\xi^L}$. 
Then the family $\{V_{\xi_i}^* e_N V_{\xi_j}\}_{i,j\in I_\xi}$ 
is a system of matrix units of $A_\xi$. 
We prepare the following projections in $N'\cap M_1$, 
\[
(z_L)_\xi :=\sum_{i\in I_\xi^L} V_{\xi_i}^* e_N V_{\xi_i}
\ \mbox{for all}\ \xi\in\Xi_L,
\quad
z_L:=\sum_{\xi\in\Xi_L}(z_L)_\xi. 
\]

In the following lemma, we determine the subfactor $L\subs M$
at the GNS Hilbert space level. 
Our argument is essentially same as the one given in \cite[Theorem 3.9]{ILP}, 
but the assumption there is different from ours. 
We present a proof for the completeness of our discussion. 

\begin{lem}
With the above settings, 
$z_L H_\vph=\ovl{\La_\vph(L)}$ holds. 
In particular, one has $z_L \in L' \cap M_1$ and $e_N\leq z_L$. 
\end{lem}
\begin{proof}
First we note that the following holds: 
\begin{equation}\label{eq: zL}
z_L H_\vph
=\ovl{\spa\{\La_\vph(V_{\xi_i}^*N)\mid \xi\in\Xi_L, i\in I_\xi^L\}}. 
\end{equation}
Indeed, let $x\in M$. Then we have 
\begin{align*}
z_L\La_\vph(x)
=&\,
\sum_{\xi\in\Xi_L}
\sum_{i\in I_\xi^L} V_{\xi_i}^* e_N V_{\xi_i}\La_\vph(x)
\\
=&\,
\sum_{\xi\in\Xi_L}
\sum_{i\in I_\xi^L} \La_\vph(V_{\xi_i}^* E_N^M(V_{\xi_i} x)). 
\end{align*}
Hence the left hand side is contained in the right hand one of (\ref{eq: zL}). 
The converse inclusion follows from 
$E_N^M(V_{\xi_i}V_{\eta_j}^*)=\de_{\xi\eta}\de_{ij}$
for $\xi,\eta\in\Xi$ and $i\in I_\xi, j\in I_\eta$. 

In particular, this yields $z_L H_\vph\subs \ovl{\La_\vph(L)}$. 
We will prove the equality 
by using the averaging technique presented 
in the proof of \cite[Thorem 3.9]{ILP} as shown below. 
To prove it, 
we may and do assume that $N, L$ and $M$ are factors 
of type III by tensoring with a type III factor. 
Assume that there would exist $x\in L$ such that $\La_\vph(x)\nin z_L H_\vph$. 
By the following equality:
\[
(1-z_L)\La_\vph(x)
=
\sum_{\xi\in\Xi\setm \Xi_L}\sum_{i\in I_\xi} 
\La_\vph(V_{\xi_i}^* E_N^M(V_{\xi_i}x))
+
\sum_{\xi\in\Xi_L}\sum_{i\in I_\xi\setm I_\xi^L} 
\La_\vph(V_{\xi_i}^* E_N^M(V_{\xi_i}x)), 
\]
the following two cases could occur: 
(I) there exists $\xi\in\Xi\setm \Xi_L$ 
such that $E_N^M(V_{\xi_i}x)\neq0$ for some $i\in I_\xi$ 
or 
(II) there exists $\xi\in \Xi_L$ 
such that $E_N^M(V_{\xi_i}x)\neq0$ for some $i\in I_\xi\setm I_\xi^L$.
In case (I), we set $I_\xi^L=\emptyset$ 
and then proceed as with case (II).
Assume that case (II) would occur. 
Take $\xi\in \Xi$ and $i\in I_\xi\setm I_\xi^L$ 
such that $E_N^M(V_{\xi_i}x)\neq0$.
Let $E_\xi\col N\ra \rho_\xi(N)$ be
the faithful normal conditional expectation with respect to $\rho_\xi$. 
By using the equality $E_N^M(V_{\xi_i}axb)=\rho_\xi(a)E_N^M(V_{\xi_i}x)b$ 
for $a,b\in N$, 
we may assume that $E_\xi(E_N^M(V_{\xi_i}x))=1$ since $N$ is of type III. 

We take a hyperfinite subfactor $R\subs N$ which is simple 
in the sense of \cite{L3}.
Then consider the weakly closed convex set 
\[
C:=\ovl{{\rm co}}^{{\rm w}}\{ux\rho_\xi(u^*)\mid u\in U(R)\}\subs L, 
\]
where $U(R)$ is the set of all unitaries in $R$. 
The hyperfiniteness of $R$ assures that 
there exists a point $w^*\in C$ such that 
$w$ satisfies $wx=\rho_\xi(x)w$ for all 
$x\in R$ and hence for all $x\in N$ by \cite[Proposition 2.10]{ILP}. 
This shows $w\in L\cap \meH_\xi=\meK_\xi$. 
Since $i\in I_\xi\setm I_\xi^L$, $V_{\xi_i}$ is orthogonal to $\meK_\xi$, 
that is, $E_N^M(V_{\xi_i}w^*)=d(\xi)^{-1}(V_{\xi_i},w)=0$. 
However $E_\xi(E_N^M(V_{\xi_i} C))=\{1\}$, and this is a contradiction. 
Therefore the cases (I) and (II) never occur, and for any $x\in L$, 
$(1-z_L)\La_\vph(x)=0$. 
This implies that $\ovl{\La_\vph(L)}\subs z_L H_\vph$. 
\end{proof}

By the previous lemma, we can describe the corner subalgebra 
$z_L M_1 z_L$ in $M_1$. 

\begin{lem}\label{lem: zM1z}
One has $z_L M_1 z_L=\ovl{L e_N L}^{\rm \,w}=Lz_L\vee \{e_N\}''$. 
\end{lem}
\begin{proof}
Recall that $e_N M_1 e_N=N e_N$. 
For $\xi,\eta\in\Xi$, $i\in I_\xi^L$ and $j\in I_\eta^L$, 
we have 
\[
V_{\xi_i}^* e_N V_{\xi_i} M_1 V_{\eta_j}^* e_N V_{\eta_j}
\subs 
V_{\xi_i}^* e_N M_1 e_N V_{\eta_j}
\subs V_{\xi_i}^* N e_N V_{\eta_j}\subs Le_N L. 
\]
This implies that $z_L M_1 z_L\subs \ovl{Le_N L}^{\, \rm w}$. 
By the previous lemma, $z_L\in L'\cap M_1$ and $z_L e_N=e_N$. 
Since $M_1$ contains $L$ and $e_N$, we have 
$z_L M_1 z_L\supset  Lz_L\vee \{e_N\}''\supset L e_N L$. 
Hence we have 
$z_L M_1 z_L=\ovl{Le_N L}^{\, \rm w}=Lz_L\vee \{e_N\}''$. 
\end{proof}

Next we will show that 
the two-step inclusion $Nz_L\subs L z_L\subs z_L M_1 z_L$
is identified with the basic extension of $N\subs L$. 
One might be able to prove this by using the abstract characterization
of the basic extension \cite[Lemma 2.4]{ILP}. 
To apply that result, we need to show that the restriction 
$\wdh{E}_M^{M_1}$ on $z_L M_1 z_L$ is an operator valued weight 
from $z_L M_1 z_L$ to $L z_L\cong L$, 
but we do not have a proof for such a statement yet. 
We avoid using this method and directly 
compare the basic extension of $N\subs L$ 
with $Nz_L\subs L z_L\subs z_L M_1 z_L$ instead. 

We set $\ps:=\om\circ E_N^L\in L_*$. 
Then $\vph|_L=\ps|_L$ holds trivially. 
Let $\{H_\ps,\La_\ps\}$ 
be the GNS representation of $L$ associated with 
the state $\ps$. 
Let $f_N\in B(H_\ps)$ be the Jones projection 
defined by 
$f_N \La_\ps(x)=\La_\ps(E_N^L(x))$ for $x\in L$. 
We set $L_1:=L\vee \{f_N\}''\subs B(H_\ps)$. 
Then we obtain the Jones' basic extension $N\subs L\subs L_1$ 
associated with $E_N^L$. 
The dual operator valued weight of $E_N^L$ is 
denoted by $\wdh{E}_L^{L_1}$. 
Note that we do not know whether $\wdh{E}_L^{L_1}$ is semifinite on 
$N'\cap L_1$ or not. 
Set a weight $\ps_1:=\ps\circ \wdh{E}_L^{L_1}\in \mP(L_1,\C)$. 
Let $\{H_{\ps_1},\La_{\ps_1}\}$ be the GNS representation of $L_1$ 
associated with the weight $\ps_1$. 
Recall the weight $\vph_1=\vph\circ \wdh{E}_M^{M_1}$ on $M_1$. 
Let $\{H_{\vph_1},\La_{\vph_1}\}$ be the GNS representation of $M_1$ 
associated with the weight $\vph_1$.
Then the following holds (\cite[Lemma 2.1]{ILP}): 
\[
H_{\ps_1}=\ovl{\La_{\ps_1}(L f_N L)},\quad 
H_{\vph_1}=\ovl{\La_{\vph_1}(M e_N M)}. 
\]
We introduce an isometry $U\col H_{\psi_1}\ra H_{\vph_1}$ satisfying 
\[
U\La_{\ps_1}(x f_N y)=\La_{\vph_1}(x e_N y), 
\quad\mbox{for}\ x, y\in L. 
\]
The well-definedness is verified as follows. 
For $x, y, a, b\in L$, we have 
\begin{align*}
\langle \La_{\vph_1}(x e_N y), \La_{\vph_1}(a e_N b) \rangle
=&\, 
\vph_1(b^* e_N a^* x e_N y)
=\vph_1(b^* E_N^M(a^* x)e_N y)
\\
=&\,
\vph\circ \wdh{E}_M^{M_1}(b^* E_N^M(a^* x)e_N y )
=
\vph(b^* E_N^M(a^* x)y)
\\
=&\,
\ps(b^* E_N^L(a^* x)y)
=
\langle 
\La_{\ps_1}(x f_N y) , \La_{\ps_1}(a f_N b)
\rangle. 
\end{align*}

\begin{lem}\label{lem: Ux}
One has 
$x U=Ux$ for $x\in L$ and 
$e_N U=U f_N$. 
\end{lem}
\begin{proof}
Since the subspace $\La_{\ps_1}(L f_N L)\subs H_{\ps_1}$ is dense, 
it suffices to show the equalities on $\La_{\ps_1}(L f_N L)$. 
Let $x,a,b\in L$. 
Then we have 
\[
x U \La_{\ps_1}(a f_N b)
=x \La_{\vph_1}(a e_N b)=\La_{\vph_1}(xa e_N b)
=U \La_{\ps_1}(xa f_N b)
=U x \La_{\ps_1}(a f_N b). 
\]
Hence $xU=Ux$. 
Next $e_N U=U f_N$ is verified as follows. 
\begin{align*}
e_N U \La_{\ps_1}(a f_N b)
=&\, 
e_N \La_{\vph_1}(a e_N b)=\La_{\vph_1}(e_N a e_N b)
\\
=&\, \La_{\vph_1}(E_N^M(a)e_N b)
=\La_{\vph_1}(E_N^L(a)e_N b)
\\
=&\, 
U \La_{\ps_1}(E_N^L(a)f_N b)
=Uf_N \La_{\ps_1}(a f_N b). 
\end{align*}
\end{proof}

Set the range projection $p_L:=UU^*\in B(H_{\vph_1})$.
It is clear that $p_L H_{\vph_1}=\ovl{\La_{\vph_1}(Le_N L)}$.
By the previous lemma or the definition of $p_L$,
$p_L$ commutes with $L$ and $e_N$.
In particular,
$p_L\in (Lz_L)' \cap \{e_N\}' \subs B(H_{\vph_1})$.

The subspace $p_L H_{\vph_1}$ plays a similar role
to the GNS Hilbert space
of $z_L M_1 z_L$ associated with the restricted weight 
$\vph_1|_{z_L M_1 z_L}$,
but $p_L H_{\vph_1}$ may not coincide
with the closure of
$\La_{\vph_1}(n_{\vph_1}\cap z_L M_1 z_L)$
because the function $t\in\R \mapsto
\si_{t}^{E_N^M \circ \wdh{E}_M^{M_1}}(z_L)\in N'\cap M_1$
may not extend to the bounded analytic function
on the strip $\{z\in\C\mid 0\leq {\rm Im}(z)\leq 1/2\}$.
However, the following lemma is sufficient for our purpose.

\begin{lem}\label{lem: zp}
In $B(H_{\vph_1})$, $p_L\leq z_L$ holds. 
\end{lem}
\begin{proof}
Using $z_L\in L'\cap M_1$ and $z_L e_N=e_N$, we have 
\[
z_L p_L H_{\vph_1}
=z_L \ovl{\La_{\vph_1}(Le_N L)}=\ovl{\La_{\vph_1}(Lz_L e_N L)}
=\ovl{\La_{\vph_1}(Le_N L)}=p_L H_{\vph_1}. 
\]
Hence $p_L\leq z_L$. 
\end{proof}

\begin{lem}
There exists an isomorphism 
$\Ps_L\col z_L M_1 z_L\ra L_1$ such that 
\begin{enumerate}

\item $\Ps_L(xz_L)=x$ for $x\in L$. 

\item $\Ps_L(e_N)=f_N$. 

\end{enumerate}
In particular, the inclusions 
$Nz_L\subs Lz_L\subs z_L M_1 z_L$ 
and $N\subs L\subs L_1$ are isomorphic. 
\end{lem}
\begin{proof}
We define the normal positive map $\Ps_L\col z_L M_1 z_L\ra B(H_{\ps_1})$ 
by $\Ps_L(x)=U^* xU$ for $x\in z_L M_1 z_L$. 
Since $p_L$ commutes with $Lz_L\vee \{e_N\}''=z_L M_1 z_L$ 
as is remarked after Lemma \ref{lem: Ux}, 
we see that $\Ps_L$ is multiplicative. 
By the previous lemma, we have 
\[
\Ps_L(z_L)=U^* z_L U=U^* z_L p_L U=U^* p_L U=1,
\]
that is, $\Ps_L$ is unital. 
Hence $\Ps_L$ is a unital $*$-homomorphism. 
By Lemma \ref{lem: Ux}, the range of $\Ps_L$ 
is equal to $U^* (Lz_L\vee \{e_N\}'') U=L\vee \{f_N\}''=L_1$. 
Also we have $\Ps_L(xz_L)=x$ for $x\in L$ and $\Ps_L(e_N)=f_N$. 
Since $z_L M_1 z_L$ is a factor, 
$\Ps_L$ is an isomorphism onto $L_1$. 
\end{proof}

Now we state our main result in this section. 

\begin{thm}\label{thm: L-discrete}
Let $N\subs M$ be an irreducible inclusion of discrete type. 
Let $N\subs L\subs M$ be an intermediate subfactor. 
Then one has the following: 
\begin{enumerate}
\item The inclusion $N\subs L$ is discrete. 

\item 
Suppose that $N$ is infinite. 
Let $\ga_N^M$ and $\ga_N^L$ be the canonical endomorphisms 
for $N\subs M$ and $N\subs L$, respectively. 
Then $[\ga_N^M|_N]$ contains $[\ga_N^L|_N]$ in $\Sect(N)$. 

\item 
Suppose that $N$ is infinite and $[\ga_N^L|_N]$ has 
in $\Sect(N)$ the following decomposition 
into irreducible sectors $[\rho_\xi]$, $\xi\in\Xi_L$:
\[
[\ga_N^L|_N]=\bigoplus_{\xi\in\Xi_L} m_\xi [\rho_\xi]. 
\]
Set 
$\meK_\xi=\{V\in L\mid Vx=\rho_\xi(x)V\ \mbox{for all}\ x\in N\}$ 
for $\xi\in\Xi_L$. 
Then one has $m_\xi=\dim(\meK_\xi)$ and 
$L$ is weakly spanned by $\meK_\xi^* N$, $\xi\in\Xi$. 
\end{enumerate}
\end{thm}
\begin{proof}
(1). 
We may and do assume that $N$ is infinite 
by tensoring with an infinite factor if necessary. 
For $\xi\in \Xi_L$, 
we define the matrix algebra $B_\xi\subs N'\cap M_1$ by
$B_\xi:=\meK_\xi^* e_N \meK_\xi$. 
Then it is easy to see that 
$z_L A_\xi z_L=(z_L)_\xi A_\xi (z_L)_\xi=B_\xi$. 
Hence we have 
\[
z_L(N'\cap M_1)z_L=\bigoplus_{\xi\in\Xi_L} B_\xi.
\] 
Let $\Ps_L\col z_L M_1 z_L\ra L_1$ be the isomorphism constructed 
in the previous lemma. 
Using the equalities 
$\Ps_L(z_L(N'\cap M_1)z_L)=\Ps_L((N z_L)'\cap z_L M_1 z_L)=N'\cap L_1$ 
and $\Ps_L(B_\xi)=\meK_\xi^* f_N \meK_\xi$, 
we have 
\[
N'\cap L_1=\bigoplus_{\xi\in\Xi_L}\meK_\xi^* f_N \meK_\xi. 
\]
Since $\wdh{E}_{L}^{L_1}$ is finite on each matrix algebra 
$\meK_\xi^* f_N \meK_\xi$, 
$\wdh{E}_L^{L_1}$ is semi-finite on $N'\cap L_1$. 
Therefore the inclusion $N\subs L$ is discrete. 

(2). 
Take $V\in \meK_\xi$ such that $V^*f_N V$ is a minimal projection 
in $\meK_\xi^* f_N \meK_\xi$. 
Note that $E_N^L(VV^*)=1$. 
The projection $V^*f_N V$ corresponds to an irreducible sector 
in $\Sect(N)$. 
The sector is actually equal to $[\rho_\xi]\in \Sect(N)$ as seen below. 
Set $W:=f_N V \in L_1$. 
Using $WW^*=f_NE_N^L(VV^*)f_N=f_N$ and $W^*W=V^*f_N V$, we have 
\[
f_N \rho_\xi(x)=WW^* \rho_\xi(x)=W V^* f_N \rho_\xi(x)
=WV^*\rho_\xi(x)f_N
=WxV^*f_N=WxW^*. 
\]
By \cite[Lemma 3.1]{ILP}, the minimal projection 
$V^* f_N V$ corresponds to $[\rho_\xi]$, 
and the canonical endomorphism 
$\ga_N^L\col L\ra N$ has the following decomposition in $\Sect(N)$ 
\[
[\ga_N^L|_N]
=\bigoplus_{\xi\in \Xi_L}\dim(\meK_\xi) [\rho_\xi]. 
\]
From this, we see that $[\ga_N^L|_N]$ is contained in $[\ga_N^M|_N]$ 
because each irreducible is contained in $[\ga_N^M|_N]$
and we trivially have $\dim(\meK_\xi)\leq \dim(\meH_\xi)$. 

(3). 
Apply \cite[Lemma 3.8]{ILP} to the discrete inclusion $N\subs L$. 
\end{proof}

\begin{rem}\label{rem: bimodule}
In the above theorem, we have shown $[\ga_N^M|_N]$ contains $[\ga_N^L|_N]$. 
This means ${}_N L^2(M)_N$ contains ${}_N L^2(L)_N$. 
Since $M$ contains $L$, it might sound trivial 
by regarding $L^2(M)\supset L^2(L)$ naturally. 
Indeed this inclusion implies that of the left $N$-modules, 
but it does not directly derive the inclusion of the right $N$-modules 
when $M$ is infinite 
because the right module structures are given by the different modular 
conjugations $J_M$ and $J_L$. 
If there exists a faithful normal conditional expectation 
$E_L^M\col M\ra L$, then we can regard the restriction $J_M$ on $L^2(L)$ 
as $J_L$ by using Takesaki's theorem \cite[p.309]{Ta}. 
Hence in this case, the above theorem is really trivial. 
The point is that the above theorem also holds 
for an intermediate subfactor $N\subs L\subs M$ such that 
there exist no faithful normal conditional expectations from $M$ onto $L$. 
\end{rem}

\section{Minimal actions of compact quantum groups}

In this section, 
we apply Theorem \ref{thm: L-discrete} to 
inclusions of factors coming from 
minimal actions of compact quantum groups, 
and we obtain the Galois correspondence (Theorem \ref{thm: Galois}). 

\subsection{Compact quantum groups}
We briefly explain compact quantum groups and their actions. 
We adopt the definition of a compact quantum group 
that is introduced in \cite[Definition 2.1]{W2} as follows: 

\begin{defn}
A \textit{compact quantum group} $\bG$ is a pair 
$(C(\bG),\delta)$ that satisfies the following conditions: 

\begin{enumerate}

\item $C(\bG)$ is a separable unital $C^*$-algebra. 

\item 
The map $\delta\col C(\bG)\ra C(\bG)\oti C(\bG)$ 
is a coproduct, i.e. it is a faithful unital $*$-homomorphism satisfying 
the coassociativity condition, 
\[(\delta\otimes\id)\circ \delta
=(\id\otimes\delta)\circ\delta.
\]

\item 
The vector spaces $\delta(\CG)(\C\otimes \CG)$ 
and $\delta(\CG)(\CG\otimes\C)$ 
are dense in $\CG\otimes \CG$.
\end{enumerate}
\end{defn}

Let $h$ be
the \textit{Haar state} on $\CG$, which satisfies the invariance 
condition: 
\[
(\id\oti h)(\de(a))=h(a)1=(h\oti\id)(\de(a))
\quad
\mbox{for all}\ a\in \CG. 
\]
In this paper, we always assume that the Haar state is faithful.
If the Haar state is tracial, we say that the compact quantum group
is of \emph{Kac type} \cite{ES}.

Let $\{\ltG,\La_{h}\}$ be the GNS representation of $\CG$ associated with $h$.
We define the von Neumann algebra $\lG=\ovl{\CG}^{\,\rm{w}}\subs B(\ltG)$.
We can extend the coproduct $\de$ to $\lG$ by the standard procedure
\cite{KV}.
The extended coproduct is also denoted by $\de$.
Then the pair $(\lG,\de)$ is a
\emph{von Neumann algebraic compact quantum group}
in the sense of \cite{KV}.

Let $H$ be a Hilbert space. 
We say that a unitary $v\in \lG\oti B(H)$
is a (left) \emph{unitary representation} on $H$ 
when it satisfies $(\de\oti\id)(v)=v_{13}v_{23}$.
The unitary representation $v$ is said to be 
\emph{irreducible} 
if $\{T\in B(H)\mid (1\oti T)v=v(1\oti T)\}=\C$. 
The set of the equivalence classes of 
all irreducible unitary representations is denoted by $\IG$. 
For $\pi\in\IG$, take a representative $v_\pi\in \lG\oti B(K_\pi)$. 
Then it is well-known that $K_\pi$ is finite dimensional. 
Set $d_\pi:=\dim(K_\pi)$. 
We denote by $\lG_\pi$ the subspace of $\lG$ 
that is spanned by the entries of $v_\pi$. 
Then the subspace $A(\bG)=\spa\{\lG_\pi\mid\pi\in\Irr(\bG)\}$ 
is a weakly dense unital $*$-subalgebra of $\lG$. 

We also use the modular objects of $\lG$.
Let $\{f_z\}_{z\in\C}$ be the Woronowicz characters
on $\AG$. 
For its characterization, readers are referred 
to \cite[Theorem2.4]{W2}. 
Then the modular automorphism group $\{\si_t^{h}\}_{t\in\R}$ on $\AG$
is given by 
\[
\si_t^{h}(x)
=(f_{it}\oti\id\oti f_{it})\big{(}(\de\oti\id)(\de(x))\big{)}
\quad\mbox{for all}\ 
t\in \R,\ x\in \AG. 
\]
We define the following map $\ta_t\col \AG\ra\AG$ by
\[
\ta_t(x)
=(f_{it}\oti\id\oti f_{-it})
\big{(}(\de\oti\id)(\de(x))\big{)}
\quad\mbox{for all}\ 
t\in \R,\ x\in \AG. 
\]
Then $\{\ta_t\}_{t\in\R}$ is a one-parameter automorphism group on $\AG$ 
and it is called the \emph{scaling automorphism group}.
Since the Haar state $h$ is invariant
under 
the $*$-preserving maps $\si_t^{h}$ and $\ta_t$, 
we can extend them to the maps on $\CG$, and on $\lG$.
By definition, we have 
\begin{equation}\label{eq: modular-scaling}
\si_t^h(x)=(f_{2it}\oti\ta_{-t})(\de(x))
\quad\mbox{for all}\ x\in \AG. 
\end{equation}

We recall the notion of a left coideal von Neumann algebra introduced 
in \cite[Definition 4.1]{ILP}. 
In this paper, we simply call it a left coideal. 

\begin{defn}
Let $B\subs \lG$ be a von Neumann subalgebra. 
We say that $B$ is a \emph{left coideal} if $\de(B)\subs \lG\oti B$.
\end{defn}

If $\bG$ comes from an ordinary compact group, 
it is known that any left coideal is of the form $L^\infty(\bG/\bH)$ 
for a closed subgroup $\bH\subs \bG$ \cite{AHKT}.
Therefore a left coideal can be considered 
as an object like a ``closed subgroup'' of $\bG$. 
Even in quantum case, 
we can also introduce the notion of a closed quantum subgroup 
$\bH\subs \bG$ as \cite{Po2}. 
However when $\bG$ is not a compact group, 
not all the left coideals of $\bG$ 
have quotient forms as $L^\infty(\bG/\bH)$ \cite{Po1,Po2,T1}. 
For a compact quantum group satisfying a certain condition, 
we have a necessary and sufficient condition 
so that a left coideal is of the form 
$L^\infty(\bG/\bH)$ \cite[Theorem 3.18]{T2}. 

Now let $B$ be a left coideal, and we put $B_\pi:=B\cap \lG_\pi$. 
Since $B$ admits the left $\bG$-action $\de$, $B$ is weakly spanned 
by $B_\pi$, $\pi\in\IG$. 

\begin{lem}\label{lem: upi}
Let $B\subs \lG$ be a left coideal and $\pi\in\IG$. 
Then there exist a unitary representation 
$u_\pi=(u_{\pi_{i,j}})_{i,j\in I_\pi}$ 
and a subset $I_\pi^B\subs I_\pi$ such that
\begin{enumerate}

\item $u_\pi$ is equivalent to $v_\pi$, 

\item $B_\pi=\spa\{u_{\pi_{i,j}}\mid i\in I_\pi,\ j\in I_\pi^B\}$. 
\end{enumerate}
\end{lem}
\begin{proof}
Let $\Hom_\bG(K_\pi, \lG)$ be the set of $\bG$-linear maps 
from $K_\pi$ into $\lG$, that is, it consists of
linear maps $S\col K_\pi\ra \lG$ such that 
$\de\circ S=(\id\oti S)\circ v_\pi$, 
where $v_\pi$ is regarded as a map from $K_\pi$
to $\lG_\pi\oti K_\pi$. 
Similarly we define $\Hom_\bG(K_\pi,B)$, which is a subspace 
of $\Hom_\bG(K_\pi,\lG)$. 
Let $(\vep_i)_{i\in I_\pi}$ be an orthonormal basis of $K_\pi$. 
We prepare the inner product of $\Hom_\bG(K_\pi,\lG)$ defined 
by
\[
\langle S|T \rangle1:=\sum_{i\in I_\pi} T(\vep_i)^* S(\vep_i).
\] 
Then $\Hom_\bG(K_\pi,\lG)$ is a Hilbert space of dimension $d_\pi$.
We take an orthonormal basis $\{S_i\}_{i\in I_\pi}$
of $\Hom_\bG(K_\pi,\lG)$ which contains an orthonormal basis
of $\Hom_\bG(K_\pi,B)$ denoted by $\{S_i\}_{i\in I_\pi^B}$.

We define the linear map $T_j\col K_\pi\ra \lG$ 
by $T_j(\vep_i)=v_{\pi_{ij}}$ for $j\in I_\pi$. 
Then it is easy to see that $\{T_j\}_{j\in I_\pi}$ is 
an orthonormal basis of $\Hom_\bG(K_\pi,\lG)$. 
Hence there exists a unitary matrix $\nu_\pi:=(\nu_{\pi_{ij}})_{i,j\in I_\pi}$ 
in $B(\C^{|I_\pi|})$ such that for $i\in I_\pi$,
\[
S_i=\sum_{j\in I_\pi} \nu_{\pi_{ji}}T_j. 
\]
We define the unitary representation 
$u_\pi:=(1\oti \nu_\pi^*)v_\pi (1\oti \nu_\pi)$. 
Then we have 
\begin{align*}
u_{\pi_{ij}}
=&\,
\sum_{k,\el\in I_\pi}(\nu_\pi^*)_{ik}v_{\pi_{k\el}} \nu_{\pi_{\el j}}
=
\sum_{k,\el\in I_\pi}(\nu_\pi^*)_{ik} \nu_{\pi_{\el j}} T_\el (\vep_k)
\\
=&\,
\sum_{k\in I_\pi}(\nu_\pi^*)_{ik} S_j (\vep_k)
=
S_j\big{(}\sum_{k\in I_\pi}(\nu_\pi^*)_{ik} \vep_k\big{)}, 
\end{align*}
and 
\[
S_j(\vep_k)=\sum_{i\in I_\pi}\nu_{\pi_{ki}}u_{\pi_{ij}}. 
\]

Therefore $u_{\pi_{ij}}\in B$ for all $i\in I_\pi$ and $j\in I_\pi^B$, 
and they span $B_\pi$. 
\end{proof}

\subsection{Minimal actions}

Let $M$ be a von Neumann algebra and $\bG$ a compact quantum group. 
Let $\al\col M\ra M\oti\lG$ be a unital faithful normal $*$-homomorphism. 
We say that $\al$ is an \emph{action} of $\bG$ on $M$ 
if $(\al\oti\id)\circ\al=(\id\oti\de)\circ\al$. 
When the subspace $\{(\ph\oti\id)(\al(M))\mid\ph\in M_*\}$ is 
weakly dense in $\lG$, we say that $\al$ has \emph{full spectrum}. 
Set $M^\al:=\{x\in M\mid \al(x)=x\oti1\}$. 
We recall the notion of minimality of an action which is introduced in 
\cite[Definition 4.3 (i)]{ILP}. 

\begin{defn}
Let $\al\col M\ra M\oti\lG$ be an action. 
We say that $\al$ is \emph{minimal} 
if $\al$ has full spectrum and satisfies $(M^\al)'\cap M=\C$. 
\end{defn}

Let $\al$ be a minimal action of $\bG$ on $M$. 
We set $N=M^\al$. 
Then the action $\al$ is dual 
when $N$ is infinite \cite[Proposition 6.4]{V1}, 
that is, there exists a $\bG$-equivariant embedding of $\lG$ into $M$. 
We can prove this result in the same line as the proof 
of \cite[Theorem 2.2]{Y1}, 
where minimal actions of compact Kac algebras have been considered. 
In particular, \cite[Proposition 2.14]{Y1} also holds 
for minimal actions of compact quantum groups. 
Hence for any $\pi\in\IG$, there exists a Hilbert space $\meH_\pi\subs M$ 
with support $1$ 
such that $\al(\meH_\pi)\subs \meH_\pi\oti \lG_\pi$. 
If $\{V_{\pi_i}\}_{i\in I_\pi}$ is an orthonormal basis of $\meH_\pi$, 
there exists $w_{\pi_{ij}}\in \lG_\pi$ for $i,j\in I_\pi$ 
such that $(V_{\pi_j}^*\oti1)\al(V_{\pi_i})=1\oti w_{\pi_{ji}}$. 
Then we see that the matrix $(w_{\pi_{i,j}})_{i,j\in I_\pi}$
is a unitary representation equivalent to $v_\pi$.
Hence we may assume that $w_{\pi_{ij}}=v_{\pi_{ij}}$ by taking the new
$(V_{\pi_i})_{i\in I_\pi}$ if necessary, that is, we have
\begin{equation}\label{eq: V-pi}
\al(V_{\pi_i})=\sum_{j\in I_\pi}V_{\pi_j}\oti v_{\pi_{ji}}.
\end{equation}
Let $\rho_{\meH_\pi}\in\End(M)$ be the endomorphism implemented 
by the Hilbert space $\meH_\pi$, that is, 
\[
\rho_{\meH_\pi}(x)=\sum_{i\in I_\pi}V_{\pi_i} x V_{\pi_i}^*
\quad\mbox{for}\ x\in M. 
\] 
Then it is easy to see that $\rho_{\meH_\pi}(N)\subs N$, 
and we denote $\rho_{\meH_\pi}|_N$ by $\rho_\pi$,
which is irreducible. 
Note that $[\rho_\pi|_{N}]\in \Sect(N)$ 
does not depend on the choice of $\meH_\pi$. 
Let $\pi,\si\in\IG$. 
Then by minimality of $\al$, 
it is shown that 
$[\rho_\pi]=[\rho_\si]\in \Sect(N)$ 
if and only if $\pi=\si$.

Define the conditional expectation $E:=(\id\oti h)\circ \al$ 
from $M$ onto $N$.
Take a faithful state $\om\in N_*$ and set $\vph:=\om\circ E\in M_*$. 
Let $\{H_\vph,\La_\vph\}$ be the GNS representation of $M$ associated 
with $\vph$.
Let $N\subs M\subs M_1:=M\vee\{e_N\}''$ be 
the basic extension where 
the Jones projection $e_N\in M_1$ is defined 
by $e_N \La_\vph(x)=\La_\vph(E(x))$ for $x\in M$. 
We denote by $\wdh{E}\in \mP(M_1,M)$ 
the dual operator valued weight associated with $E$. 

Now we assume that $N$ is infinite. 
Take $\meH_\pi$ as before. 
Set $A_\pi:=\meH_\pi^* e_N \meH_\pi$. 
Then $A_\pi$ is contained in $N'\cap M_1$, 
and $\{A_\pi\}_{\pi\in\IG}$ are $d_\pi\times d_\pi$-matrix algebras, 
respectively. 
Moreover on $A_\pi$, the weight $E\circ\wdh{E}$ is finite. 
Let $z_\pi\in A_\pi$ be the unit projection. 

\begin{lem}
When $M^\al$ is infinite, 
the following statements hold: 
\begin{enumerate}

\item $1=\displaystyle \sum_{\pi\in\IG}z_\pi$. 

\item $\displaystyle N'\cap M_1=\bigoplus_{\pi\in\IG}A_\pi$. 

\item The inclusion $N\subs M$ is discrete. 
\end{enumerate}
\end{lem}
\begin{proof}
(1). 
Take $\{W_{\pi_k}\}_{k=1}^{d_\pi}$ in $\meH_\pi$ 
such that $E(W_{\pi_k} W_{\pi_\el}^*)=\de_{k,\el}1$. 
Then we have 
\[\displaystyle z_\pi=\sum_{k=1}^{d_\pi} W_{\pi_k}^* e_N W_{\pi_k}.
\] 
Since $M$ is weakly spanned by $\meH_\pi^* N$, $\pi\in\IG$, 
we have 
\[
H_\vph=\ovl{\spa \{ \La_\vph(\meH_\pi^* N )\mid \pi\in\IG\}}.
\] 
For any $x\in N$ and $V_\si \in \meH_\si$ with $\si\neq \pi$, we have 
\[
z_\pi \La_\vph(V_\si^* x ) 
=\sum_{k=1}^{d_\pi} W_{\pi_k}^* \La_\vph(E(W_{\pi_k} V_\si^*)x)=0, 
\]
where we have used $E(\meH_\pi\meH_\si^*)=0$ in the last equality. 
If $\si=\pi$, we have $z_\pi \La_\vph(V_\pi^* x)=\La_\vph( V_\pi^* x)$. 
Hence the range space of $z_\pi$ coincides with $\La_\vph(\meH_\pi^* N)$, 
and $\{z_\pi\}_{\pi\in\IG}$ is a partition of unity.

(2).
We first show that $z_\pi$ is a central projection in $N'\cap M_1$. 
For $\pi\in\IG$, take $\{W_{\pi_k}\}_{k=1}^{d_\pi}$ as above. 
It suffices to prove that $z_\pi (N'\cap M_1)z_\si=0$ if $\pi\neq\si$. 
Let $x\in N'\cap M_1$ and take $x_0\in N$ such that 
$x_0 e_N=e_N W_{\pi_k}xW_{\si_\el}^*e_N$. 
Then for any $y\in N$, we have 
\begin{align*}
x_0\rho_\si(y)e_N
=&\,
e_N W_{\pi_k}xW_{\si_\el}^* \rho_\si(y)e_N
=
e_N W_{\pi_k}xyW_{\si_\el}^*e_N
\\
=&\,
e_N W_{\pi_k}yxW_{\si_\el}^*e_N
=
e_N \rho_\pi(y)W_{\pi_k}xW_{\si_\el}^* e_N
\\
=&\,
\rho_\pi(y)x_0e_N. 
\end{align*}
This shows that $x_0$ intertwines $\rho_\si$ and $\rho_\pi$. 
So, we get $x_0=0$. 
Hence we have 
$e_N W_{\pi_k}(N'\cap M_1)W_{\si_\el}^*e_N=0$ 
and $z_\pi (N'\cap M_1)z_\si=0$. 

Second we show that each $p_k:=W_{\pi_k}^* e_N W_{\pi_k}$ 
is a minimal projection in $N'\cap M_1$. 
This is because 
the reduced inclusion $Np_k\subs p_k M_1 p_k$ is isomorphic to 
the irreducible inclusion $\rho_\pi(N)\subs N$. 
Hence $N'\cap M_1$ is the direct sum of $A_\pi$, $\pi\in\IG$. 

(3). This is trivial by (2) because $\wdh{E}$ is finite on each $A_\pi$. 
\end{proof}

By the previous lemma, we can regard $\Xi=\IG$. 

\begin{defn}
Let $\al\col M\ra M\oti\lG$ be a minimal action. 
\begin{enumerate}
\item 
For an intermediate subfactor $M^\al\subs L \subs M$, 
we define the weakly closed subspace $\meL(L)\subs \lG$ 
by 
\[
\meL(L)=\ovl{\spa}^{{\rm w}}\{(\om\oti\id)(\al(L))\mid \om\in M_*\}. 
\]
\item 
For a left coideal $B\subs \lG$, 
we define the intermediate subfactor $M^\al\subs \meM(B)\subs M$ by 
\[
\meM(B)=\{x\in M\mid \al(x)\in M\oti B\}. 
\]
\end{enumerate}
\end{defn}

We also denote by $\meL_\al(L), \meM_\al(B)$ for $\meL(L), \meM(B)$ 
when we want to specify the action $\al$. 
From now, we freely use the notations prepared in the previous 
section. 

\begin{lem}
For any intermediate subfactor $M^\al\subs L\subs M$, 
$\meL(L)$ is a left coideal of $\bG$. 
\end{lem}
\begin{proof}
If we consider a minimal action $\be:=\id\oti\al$ on $B(\el_2)\oti M$, 
then $\meL_\be(B(\el_2)\oti L)=\meL_\al(L)$. 
Therefore we may assume that $M^\al$ is infinite. 
We have to check that
$\de(\meL(L))\subs \lG\oti \meL(L)$
and $\meL(L)$ is multiplicatively closed.

Set $\meK_\pi=\meH_\pi\cap L$. 
Then $L$ is weakly spanned by 
$\meK_\pi^* N$ by Theorem \ref{thm: L-discrete}. 
Recall two bases $\{V_{\pi_i}\}_{i\in I_\pi}$ 
and $\{W_{\pi_i}\}_{i\in I_\pi}$ in $\meH_\pi$ as before, that is,
we have the equalities 
$V_{\pi_i}^* V_{\pi_j}=\de_{ij}$, (\ref{eq: V-pi})
and $E(W_{\pi_i} W_{\pi_j}^*)=\de_{ij}$.
We may assume that $\{W_{\pi_i}\}_{i\in I_\pi^L}$ is a basis 
of $\meK_\pi$. 
There exist $c_{\pi_{j i}}\in\C$, $i,j\in I_\pi$ such that 
$V_{\pi_j}^*W_{\pi_i}=c_{\pi_{j i}} $. 
Note that the matrix $(c_{\pi_{j i}})_{i,j\in I_\pi}$ is 
invertible. 
Since 
\begin{align}
\al(W_{\pi_i})
=&\,
\sum_{j\in I_\pi}\al(V_{\pi_j}V_{\pi_j}^* W_{\pi_i})
=
\sum_{j\in I_\pi}c_{\pi_{ji}}\al(V_{\pi_j})
\notag\\
=&\,
\sum_{j,k\in I_\pi}
c_{\pi_{j i}} V_{\pi_k}\oti v_{\pi_{kj}}
=
\sum_{k\in I_\pi}
V_{\pi_k}\oti \big{(}
\sum_{j\in I_\pi}c_{\pi_{j i}} v_{\pi_{kj}}
\big{)}
, \label{eq: al-W-pi}
\end{align}
we have 
\begin{equation}\label{eq: LofL}
\meL(L)
=\ovl{\spa}^{{\rm w}}
\big{\{}
\sum_{j\in I_\pi}c_{\pi_{ji}}v_{\pi_{k,j}}
\mid i\in I_\pi^L,\ j,k\in I_\pi,\ \pi\in\Xi_L
\big{\}}. 
\end{equation}
This implies that $\de(\meL(L))\subs \lG\oti \meL(L)$. 
Let $\pi,\si\in\IG$. 
Next we show that $\meL(L)$ is multiplicatively closed. 
It suffices to show that the product of 
$\displaystyle\sum_{j\in I_\pi}c_{\pi_{ji}}v_{\pi_{rj}}$ and 
$\displaystyle\sum_{j\in I_\si}c_{\si_{j\el}}v_{\si_{sj}}$ is 
contained in $\meL(L)$ 
for all $(i,r)\in I_\pi^L\times I_\pi$ 
and $(\el,s)\in I_\si^L\times I_\si$.
Then it is clear because the left hand side of the following equality 
is contained in $\C\oti \meL(L)$:
\[
(V_{\si_s}^*V_{\pi_r}^*\oti1)\al(W_{\pi_i}W_{\si_\el})
=
1\oti 
\big{(}\sum_{j\in I_\pi}c_{\pi_{j i}}v_{\pi_{rj}}\big{)}
\big{(}\sum_{j\in I_\si}c_{\si_{j \el}}v_{\si_{sj}}\big{)}
\]
\end{proof}

We present a Galois correspondence 
which is a generalization of \cite[Theorem 4.4]{ILP} 
to minimal actions of compact quantum groups. 

\begin{thm}[Galois correspondence]\label{thm: Galois}
Let $\bG$ be a compact quantum group and $M$ a factor. 
Let $\al\col M\ra M\oti\lG$ be a minimal action. 
Then there exists an isomorphism 
between the lattice of intermediate subfactors of $M^\al\subs M$ 
and the lattice of left coideals of $\bG$. 
More precisely, the maps $\meM$ and $\meL$
are the mutually inverse maps, that is, 
for any intermediate subfactor $M^\al\subs L\subs M$ and 
any left coideal $B\subs \lG$, one has 
\[
\meM(\meL(L))=L,\quad \meL(\meM(B))=B. 
\]
\end{thm}
\begin{proof}
If we consider a minimal action $\be:=\id\oti \al$ on $B(\el_2)\oti M$, 
then we have $\meL_\be(B(\el_2)\oti L)=\meL_\al(L)$ 
and $\meM_\be(B)=B(\el_2)\oti \meM_\al(B)$. 
Hence we may and do assume that $M^\al$ is infinite. 

By definition, we see that $L\subs \meM(\meL(L))$. 
We will show $\meM(\meL(L))\subs L$. 
Set $\meK_\pi:=\meH_\pi\cap L$ 
and $\wdt{\meK}_\pi:=\meH_\pi\cap \meM(\meL(L))$. 
Then $\meM(\meL(L))$ is $\si$-weakly spanned by $M^\al \wdt{\meK}_\pi$ 
for $\pi\in\IG$ by Theorem \ref{thm: L-discrete}. 
We choose a basis  
$\{W_{\pi_i}\}_{i\in I_\pi}$ in $\meH_\pi$ such that 
$E(W_{\pi_i}W_{\pi_j}^*)=\de_{i,j}$ as before. 
We may assume that it contains bases of $\meK_\pi$ and $\wdt{\meK}_\pi$, 
which are denoted by $\{W_{\pi_i}\}_{i\in I_\pi^L}$ and 
$\{W_{\pi_i}\}_{i\in J_\pi}$, respectively. 
We use the invertible matrix $(c_{\pi_{i,j}})_{i,j\in I_\pi}$ 
as in the previous lemma. 

Let $j\in J_\pi$. 
Since $W_j\in \meM(\meL(L))$, 
$\al(W_j)$ is contained in $M\oti \meL(L)$, that is,
$\displaystyle\sum_{k\in I_\pi}c_{\pi_{kj}}v_{\pi_{\el k}}\in\meL(L)$ 
for all $\el\in I_\pi$ by (\ref{eq: al-W-pi}).
Recall that the $(v_{\pi_{k\el}})_{k,\el\in I_\pi}$ are linearly independent. 
By (\ref{eq: LofL}), 
there exists $d_{\pi_{ij}}\in\C$ for $i\in I_\pi^L$ 
such that for any $\el\in I_\pi$, 
\[
\sum_{k\in I_\pi}c_{\pi_{kj}}v_{\pi_{\el k}}
=
\sum_{i\in I_\pi^L}d_{\pi_{ij}}
\left(\sum_{k\in I_\pi} c_{\pi_{ki}}v_{\pi_{\el k}}\right), 
\]
that is, 
\begin{equation}\label{eq: c-cd}
c_{\pi_{kj}}=\sum_{i\in I_\pi^L}c_{\pi_{ki}}d_{\pi_{ij}}
\quad\mbox{for all}\ j\in J_\pi,\ k\in I_\pi.
\end{equation}
Note that $d_{\pi_{ij}}$ does not depend on $\el$. 
We know the matrix $C:=(c_{\pi_{k \el}})_{k\el\in I_\pi}$ is invertible. 
Multiplying $(C^{-1})_{\el k}$ ($\el\in I_\pi^L$) 
to the both sides of the above equality, summing up with $k$, we have 
\[
\de_{\el j}=d_{\pi_{\el j}}\quad\mbox{for all}\ \el\in I_\pi^L. 
\]
This yields $j\in I_\pi^L$. 
Indeed, if $j\nin I_\pi^L$, then $d_{\pi_{\el j}}=0$ for all $\el\in I_\pi^L$. 
Together with (\ref{eq: c-cd}), we have $c_{\pi_{kj}}=0$ for all $k\in I_\pi$. 
Then we have 
\[
W_{\pi_j}
=\sum_{k\in I_\pi}V_{\pi_k}(V_{\pi_{k}}^* W_{\pi_{j}})
=\sum_{k\in I_\pi}V_{\pi_k}c_{\pi_{kj}}
=0,
\]
but this is a contradiction.
Therefore $W_{\pi_j}\in L$ for any $j\in J_\pi$, and $\meM(\meL(L))\subs L$.

Next we will show that $\meL(\meM(B))=B$.
By definition, the inclusion $\meL(\meM(B))\subs B$ holds. 
We prove $B\subs \meL(\meM(B))$.
Since $B$ is $\si$-weakly spanned by subspaces 
$B_\pi=B\cap \lG_\pi$, $\pi\in\IG$,
it suffices to show that $B_\pi\subs \meL(\meM(B))$ for any $\pi\in\IG$. 
By Lemma \ref{lem: upi},
there exists a unitary matrix 
$\nu_\pi=(\nu_{\pi_{ij}})_{ij\in I_\pi}\in B(\C^{|I_\pi|})$ 
such that
$B_\pi$ is spanned by $u_{\pi_{ij}}$, $i\in I_\pi$ and $j\in I_\pi^B$, 
where $u_\pi=(1\oti \nu_\pi^*)v_\pi (1\oti \nu_\pi)$. 
For $i\in I_\pi$, we put 
$\displaystyle V_{\pi_i}':=\sum_{j\in I_\pi}\nu_{\pi_{ji}}V_{\pi_j}$. 
Then we have 
\begin{align*}
\al(V_{\pi_i}')
=&\,
\sum_{j\in I_\pi}\nu_{\pi_{ji}}\al(V_{\pi_j})
=\sum_{j,k\in I_\pi}\nu_{\pi_{ji}}(V_{\pi_k}\oti v_{\pi_{kj}})
\\
=&\,
\sum_{k\in I_\pi}(V_{\pi_k}\oti (v_\pi(1\oti\nu_\pi))_{ki})
=
\sum_{k\in I_\pi}(V_{\pi_k}\oti ((1\oti\nu_\pi)u_\pi)_{ki})
\\
=&\,
\sum_{j, k\in I_\pi}(\nu_{\pi_{kj}}V_{\pi_k}\oti u_{\pi_{ji}})
=
\sum_{j\in I_\pi} V_{\pi_j}'\oti u_{\pi_{ji}}. 
\end{align*}

Let $i\in I_\pi^B$. 
Then $u_{\pi_{ji}}\in B_\pi$ for all $j\in I_\pi$,
and $V_{\pi_i}'\in \meM(B)$ 
by the above equality. 
Again by the above equality, $u_{\pi_{ji}}\in \meL(\meM(B))$
for all $j\in I_\pi$. 
This implies $B_\pi\subs \meL(\meM(B))$. 
\end{proof}

When $\bG$ is of Kac type, it has been proved in \cite{ILP} 
that there exists a normal conditional expectation 
from $M$ onto any intermediate subfactor of $M^\al\subs M$. 
However, if $\bG$ is not of Kac type, then this is not the case 
in general as we will see below. 
We can characterize which intermediate subfactor has 
such a property. 
We recall the following notion introduced in
\cite[Definition 3.1 (2)]{T2}

\begin{defn}
Let $B\subs \lG$ be a left coideal. 
We say that $B$ has the \emph{expectation property} 
if 
there exists a faithful normal conditional expectation 
$E_B$ from $\lG$ onto $B$ satisfying $h\circ E_B=h$. 
\end{defn}

The following lemma is probably well-known for specialists. 
Since we could not find it in the literature, 
we present a proof. 
\begin{lem}\label{lem: exp-property}
Let $B\subs \lG$ be a left coideal. 
Then the following statements are equivalent:
\begin{enumerate}

\item $B$ has the expectation property. 

\item $\si_t^h(B)=B$ for all $t\in\R$.

\item $\ta_t(B)=B$ for all $t\in\R$. 
\end{enumerate}
\end{lem}
\begin{proof}
(1)$\Rightarrow$(2). 
This follows from Takesaki's theorem \cite[p.309]{Ta}.

(2)$\Rightarrow$(3). 
Since $B_\pi$, $\pi\in\IG$ spans a dense subspace of $B$, 
it suffices to show that $\ta_t(B_\pi)\subs B_\pi$ 
for all $t\in\R$ and $\pi\in\IG$. 
Recall the equality (\ref{eq: modular-scaling}). 
Then for $x\in B_\pi$, we have 
$\ta_t(x)=(f_{2it}\oti\si_{-t}^h)(\de(x))$. 
Since $B$ is a left coideal globally invariant
under the modular group $\si^h$,
we see that $\ta_t(x)\in B_\pi$.
Hence $\ta_t(B_\pi)\subs B_\pi$ for all $t\in\R$ and $\pi\in\IG$.

(3)$\Rightarrow$(1).
Let $\pi\in\IG$ and $x\in B_\pi$.
By (\ref{eq: modular-scaling}), we have
$\si_t^h(x)=(f_{2it}\oti\ta_{-t})(\de(x))$. 
Since $B$ is a left coideal globally invariant under
the scaling group $\ta$,
we see that $\si_t^h(x)\in B_\pi$.
Hence $\si_t^h(B_\pi)\subs B_\pi$ for all $t\in\R$ and $\pi\in\IG$, 
and $B$ is globally invariant under the modular group $\si^h$.
Again by Takesaki's theorem, there exists a faithful normal
conditional expectation $E_B\col \lG\ra B$ preserving $h$.
Hence $B$ has the expectation property.
\end{proof}

\begin{thm}\label{thm: exp-prop}
Let $\al$ be a minimal action of $\bG$ on a factor $M$. 
Let $M^\al\subs L\subs M$ be an intermediate subfactor. 
Then 
there exists a faithful normal conditional expectation 
$E_L^M\col M\ra L$ 
if and only if the left coideal $\meL(L)$ has the expectation property. 
\end{thm}
\begin{proof}
Let $\om\in N_*$ be a faithful state. 
Put $\vph:=\om\circ E\in M_*$. 
We note that 
$L$ is the image of a faithful normal conditional expectation 
of $M$ if and only if 
$\si_t^\vph(L)\subs L$ for all $t\in\R$. 
Indeed, if the former condition holds, 
there exists a faithful normal conditional expectation 
$E_L^M\col M\ra L$. 
Then the conditional expectation 
$E_N^L\circ E_L^M$ is equal to $E$ because $N\subs M$ is 
irreducible. 
Hence $\vph\circ E_L^M=(\vph\circ E_N^L)\circ E_L^M=\vph\circ E=\vph$. 
Then by Takesaki's theorem 
\cite[p.309]{Ta}, the latter condition holds. 
The converse implication also follows from his theorem. 

Since $\vph$ is invariant under the action $\al$, 
we have $\al\circ \si_t^\vph=(\si_t^\vph\oti\ta_{-t})\circ\al$ 
for all $t\in\R$ by \cite[Th\'{e}or\`{e}me 2.9]{E}. 
Put $B:=\meL(L)$. 
Then we have 
$L=\{x\in M\mid \al(x)\in M\oti B\}$ by Theorem \ref{thm: Galois}. 
So, for $x\in L$, 
$\si_t^\vph(x)\in L$ if and only if $\al(x)\in M\oti \ta_t(B)$. 
The von Neumann subalgebra $\ta_t(B)$
is also a left coideal by the equality
$\de\circ\ta_t=(\ta_t\oti\ta_t)\circ \de$.
Hence $\si_t^\vph(x)\in L$ if and only if $x\in \meM(\ta_t(B))$. 
Therefore, 
$\si_t^\vph(L)\subs L$ if and only if $L\subs \meM(\ta_t(B))$. 
Since $L=\meM(B)$, this is equivalent with $B\subs \ta_t(B)$. 
Hence $L$ is the image of a faithful normal conditional expectation 
of $M$ iff $B=\ta_t(B)$ for all $t\in\R$.
By the previous lemma, this equivalently means that
$B$ has the expectation property. 
\end{proof}

If $\bG$ is of Kac type, the Haar state $h$ is a faithful trace. 
Hence any left coideal has the expectation property. 
Then we have the following result which has been already shown 
in \cite[Theorem 3.9]{ILP}. 

\begin{cor}
Let $\al$ be a minimal action of $\bG$ on a factor $M$. 
Let $M^\al\subs L\subs M$ be an intermediate subfactor. 
If $\bG$ is of Kac type, then there exists a faithful normal 
conditional expectation from $M$ onto $L$. 
\end{cor}

\begin{ex}
We consider the twisted $SU(2)$ group, $\suq$ \cite{W1} and 
its minimal action $\al$ on a full factor $M$
as constructed by Ueda \cite{U1}.
Then by minimality of $\al$, intermediate subfactors 
bijectively correspond to left coideals.
By using Lemma \ref{lem: exp-property}, we can show 
the quantum spheres $L^\infty(S_{q,\th}^2)$ with $0<\th\leq \pi/2$ 
\cite{Po1,Po2} are left coideals 
without expectation property. 
By Theorem \ref{thm: exp-prop}, 
there are no faithful normal conditional expectations 
from $M$ onto the corresponding subfactors.
\end{ex}

\noindent\textbf{Acknowledgments.} 
The author would like to thank Masaki Izumi for informing him of 
Remark \ref{rem: bimodule}. 
The author is also grateful to Yoshimichi Ueda 
for various useful comments.

\end{document}